\documentclass[12pt,leqno]{article}
\def\s{\dot{s}} 
\def\w{\dot{w}}
\def\v{{\rm v}}
\def\sigmad{\dot{\sigma}} 
\def\int{\mathbb{Z}}

\def\Ue{{\cal U}_{\varepsilon}({\mathfrak g})}

\def\proof{{\bf Proof. }}

\def\pf{\proof}

\def\V{{\cal V}}
\def\m{\dot{w}_C}
\usepackage{times}
\usepackage{graphics}
\usepackage{amssymb}
\usepackage{amscd}
\usepackage{amsmath}
\input xy 
\xyoption{all}
\title{On spherical twisted conjugacy classes,third version}
\newtheorem{theorem}{Theorem}[section]
\newtheorem{lemma}[theorem]{Lemma}
\newtheorem{corollary}[theorem]{Corollary}
\newtheorem{proposition}[theorem]{Proposition}
\newtheorem{definition}[theorem]{Definition}
\newtheorem{remark}[theorem]{Remark}

\author{Giovanna Carnovale\\
Dipartimento di Matematica Pura ed Applicata\\
via Trieste 63 - 35121 Padova - Italy\\
email: carnoval@math.unipd.it }

\date{}

\begin{document}
\maketitle
\begin{abstract}
Let $G$ be a simple algebraic group over an algebraically closed field of good odd characteristic, and let $\theta$ be an automorphism of $G$ arising from an involution of 
its Dynkin diagram. We show that 
the spherical $\theta$-twisted conjugacy classes are precisely those intersecting only Bruhat cells corresponding to twisted involutions in the Weyl group. We show how the analogue of this  statement fails in the triality case. As a by-product, we obtain a dimension formula for spherical twisted conjugacy classes that was originally obtained by J-H. Lu in characteristic zero. 
\end{abstract}



\section{Introduction}

Twisted conjugacy classes were originally introduced by Gantmakher in \cite{gant} and developed in \cite{VO}, where they were viewed as orbits under the  conjugacy action of the identity component in a disconnected algebraic group. It is needless to mention that reductive disconnected groups frequently occur 
in the study of algebraic groups, for example, as centralizers of semisimple elements in non-simply-connected semisimple groups. Twisted classes also occur in problems concerning conjugacy classes of rational forms.

In recent years the attention to the twisted conjugacy classes of  an algebraic group $G$ has increased in different contexts of mathematics: for example the closure of a twisted Steinberg fiber in the wonderful compactification of a simple linear algebraic group has been computed in  \cite{he-thomsen} and twisted conjugacy classes have been shown to be Poisson submanifolds with respect to a natural Poisson structure $\pi_\theta$ induced by an automorphism $\theta$ of $G$ (\cite{lu-yakimov}).   Moreover, conjugacy classes in disconnected groups play a role in physics, due to their  connection with branes in the Wess-Zumino-Witten model (see, for instance \cite{fffs}, where they are called twined conjugacy classes). In a different context, twisted conjugacy classes in finite simple groups, so also in simple groups of Lie type, occur in the classification of racks, which is an important tool for the classification of finite dimensional pointed Hopf algebras (\cite{AG}).

Given an automorphism $\theta$ of $G$, the simplest example of a  twisted conjugacy class is the class $G*1=\{g\theta(g)^{-1},\,g\in G\}$ of the unit element in $G$. When $\theta$ 
is an involution, $G*1$ has been extensively studied in \cite{R,results}. It provides a model for symmetric spaces, and it is shown in \cite{results,vust,mauro} that a Borel subgroup $B$ of $G$ acts on this class with finitely many orbits. Transitive $G$-varieties satisfying this property are called {\em spherical}. 
The combinatorics of the Zariski closures of the $B$-orbits in $G*1$ has been described in \cite{RS} by means of a map from the set of $B$-orbits to the (set of twisted involutions in the)  Weyl group. This map is given by looking at which Bruhat cell contains the twisted $B$-orbit. In the untwisted case, the analysis of the intersection of conjugacy classes and Bruhat cells is a powerful tool and it has been used in many different situations. For spherical twisted conjugacy classes, this analysis is subject of current research. We summarize here the recent developments.

A characterization of spherical $\theta$-twisted conjugacy classes, when the automorphism $\theta$ is induced from an automorphism of the Dynkin diagram of $G$ and the characteristic of the base field is zero, is given in \cite{lu}. It is provided in terms of a dimension formula involving the Weyl group element whose associated  Bruhat cell intersects a class $C$ densely.  Such a Weyl group element, which we shall denote by $w_C$, is the maximum  among all Weyl group elements $\sigma$ for which $C\cap B\sigma B$ is non-empty. The above mentioned characterization generalizes to the twisted case a result described in \cite{ccc,gio-mathZ} with a more elegant proof, but it requires some restrictions on the base field. 
The motivation for such a generalization lies in the relation of the element $w_C$ with the smallest dimension of symplectic leaves of the natural Poisson structure $\pi_\theta$ on $G$. Besides, the  dimension formula is related to the vanishing of $\pi_\theta$ in a class $C$.

The aim of the present paper is to provide another characterization of spherical $\theta$-twisted conjugacy classes, when $\theta$ is induced from an involution of the Dynkin diagram, by means of their intersection with Bruhat cells. This has to be seen as a twisted analogue of some results in \cite{gio-mathZ} and the main result in \cite{gio-fourier}. It can be formulated as follows, when we restrict to simply-connected groups. 

\smallskip

\noindent{\bf Theorem}  {\em Let $G$ be a simply-connected simple algebraic group over an algebraically closed field of good odd characteristic. Let $\theta$ be an automorphism of $G$ induced by an involution of its Dynkin diagram. A $\theta$-twisted conjugacy class is spherical if and only if it intersects only Bruhat cells corresponding to twisted involutions in the Weyl group of $G$.}

\smallskip

The triality case falls out of this picture. Indeed, we show that there are no twisted classes intersecting only Bruhat cells corresponding to twisted involutions in the Weyl group, whereas it is shown in \cite{mauro, lu} that there exists a spherical twisted conjugacy class. We expect that the combination of this characterization with the one in \cite{lu} can be exploited in order to  obtain a complete classification of spherical twisted conjugacy classes when $\theta$ is an involution of the Dynkin diagram. This is part of a forthcoming project. 

As a by-product of our results, we are able to prove Lu's dimension formula when $\theta$ is an involution and $k$ is of good, odd characteristic. This can be stated, for $G$ simply-connected, as follows:

\smallskip

\noindent{\bf Theorem}  {\em Let $G$ be a simply-connected simple algebraic group over an algebraically closed field of good odd characteristic. Let $\theta$ be an automorphism of $G$ induced by an involution of its Dynkin diagram $\theta'$. A $\theta$-twisted conjugacy class $C$ is spherical if and only if  $$\dim C=\ell(w_C)+{\rm rk}(1-w_C\theta').$$}
%
%
Here $\ell$ denotes the length function in the Weyl group and ${\rm rk}$ denotes the rank of the operator in the geometric representation of the Weyl group.

It was pointed out in \cite[Remark 1.2]{lu} that the case of $\tau$-twisted classes for a general automorphism $\tau$ of $G$ can be reduced to the above setting as follows. 
For a fixed maximal torus $T$ contained in $B$ we have the equality  $(\tau(B),\tau(T))=(gBg^{-1}, gT g^{-1})$ for some $g\in G$. Multiplying $g$ on the right by a suitable element in $T$, we may choose $g$ so that 
$\theta:= {\rm Int} (g^{-1})\circ \tau$ is the automorphism of $G$ induced from an automorphism of the Dynkin diagram. Then,  right translation by $g$ induces a $G$-equivariant isomorphism between the $\tau$-twisted conjugacy class of an element $x$ and the $\theta$-twisted conjugacy class of $xg$. 

\smallskip

The paper is structured as follows. The basic notation and terminology, and the first properties of twisted conjugacy classes are provided in Sections \ref{notation} and \ref{twisted-classes}. The first properties of spherical twisted conjugacy classes are dealt with in Section \ref{spherical-classes}. Here, it is shown that if $\theta$ is an involution, then a spherical $\theta$-twisted conjugacy class intersects only Bruhat cells associated with $\theta$-twisted involutions in the Weyl group of $G$. 
This result is obtained by induction on the length of a path in the set  $\V$ of $B$-orbits which is constructed using the action on $\V$,  defined in \cite{RS}, of a monoid associated with the Weyl group,  and the Weyl group action on $\V$  introduced in \cite{knop}. The approach is similar to that in \cite{gio-mathZ} but the proof has been shortened and simplified.
In Section \ref{involutive} we analyze the twisted conjugacy classes intersecting only Bruhat cells corresponding to twisted involutions in the Weyl group (involutive classes). By a simple case-by-case analysis  on the possible maximal elements $w_C$'s we get to a better understanding of a representative lying in the Bruhat cell corresponding to $w_C$. Here, we use the classification of all possible $w_C$'s in \cite{lu}, which holds under very mild restrictions on the base field. The case-by-case analysis is simpler here than in \cite{gio-mathZ} because there are less cases to be dealt with. In Section \ref{quasi} we show that, except from the case in which $w_C=w_0$ and $G$ is of type $D_{2n}$,  if $C$ is an involutive twisted conjugacy class, then there are finitely many $B$-orbits in $B w_C B$. A simple topological argument  shows that  $C$ is spherical. The strategy is similar to the strategy used in \cite{gio-fourier} in order to deal with the case $w_C=w_0=-1$. However, it has been improved in order to be applied to a wider range of cases, namely all but the one in which $w_C=w_0$ but $w_0\neq-\theta$. The remaining case is dealt with in Section \ref{remaining}. Here we need to use a different argument. We show that,  for a suitable representative $x$  of an involutive class $C$, with stabilizer $G_x$, the set $BG_x$ is dense in $G$. We do so  by showing that the intersection of $G_x B$ with $U\sigma B$ is dense in $U\sigma B$ for every $\sigma$ in the Weyl group. This concludes the proof when $w_C=w_0$ and $G$ is of type $D_{2n}$. Here, the final strategy resembles the strategy used in \cite[Section 5]{gio-fourier}. However, the techniques used in Lemma ~\ref{piccola} are specific of the twisted case and \cite[Theorem 5.7]{gio-fourier} is extended in Lemma~\ref{dense} to a statement on a general transitive $G$-variety. On the other hand, the computational work in this paper is considerably simpler than in \cite{gio-mathZ,gio-fourier} because we do not need to consider doubly-laced root systems. Finally, in  Section \ref{dimension-section} we show how to apply the obtained results in order to retrieve Lu's dimension formula in good odd characteristic, when $\theta$ is an involution.

\section{Notation}\label{notation}

Unless otherwise stated, $G$ is a simply-connected, simple algebraic group over an algebraically closed field $k$ of zero or odd good characteristic. We recall that the characteristic  is good if it does not divide the coefficients in the expression of the highest root as a linear combination of simple roots. 
Let $T$ be a fixed maximal torus of $G$,  and let $\Phi$ be the associated
root system. Let $B\supset T$ be a Borel subgroup with unipotent radical $U$, let
$\Delta=\{\alpha_1,\ldots,\alpha_n\}$ be the basis of $\Phi$ relative
to $(T,\,B)$, with numbering of the simple roots as in \cite{bourbaki}. The set of positive roots will be denoted by $\Phi^+$. The Weyl group of $G$ will be denoted by $W$ and the reflection with respect to $\alpha\in\Phi$ will be denoted by $s_\alpha$. For $w\in W$ we shall denote by $\w$ a representative of $w$ in $N(T)$. 
The symbol $X_\alpha$ will denote the root subgroup corresponding to $\alpha$. We will choose parametrisations of roots subgroups  $x_\alpha(\xi)$ and $x_{-\alpha}(\xi)$ for $\alpha\in\Phi^+$  in such a way that 
$\s_\alpha=x_\alpha(1)x_{-\alpha}(-1)x_\alpha(1)$ lies in $N(T)$ (\cite[Lemma 8.1.4]{springer}). 
Then, as in \cite{yale}, for any root  $\alpha\in\Phi^+$ and any $\xi\in k^*$ we define $h_\alpha(\xi) = x_\alpha(\xi)x_{-\alpha}(-\xi^{-1})x_\alpha(\xi)\s_\alpha^{-1}\in T$.
By $\theta$ we denote a non-trivial automorphism of the Dynkin diagram of $G$.  
By abuse of notation, the induced automorphism of $G$ will also be denoted by $\theta$. 
We recall that for every $\alpha\in\Phi$ there is $\epsilon_\alpha\in\{\pm1\}$ such that $\theta(x_\alpha(\xi))=x_{\theta\alpha}(\epsilon_\alpha\xi)$, with $\epsilon_\alpha=1$ for $\alpha\in\Delta\cup(-\Delta)$ (see \cite[Corollary to Theorem 29]{yale}). By \cite[Lemma 8.1.4(iv)]{springer}, we deduce that $\epsilon_{\beta}=\epsilon_{-\beta}$ for every $\beta\in\Phi^+$ so $\theta(h_\beta(\xi))=h_{\theta(\beta)}(\xi)$. It was observed in \cite[Proposition 2.1]{mohr} that, unless $\Phi$ is of type $A_{2n}$, one may choose the parametrisations $x_{\pm\gamma}(\xi)$  in such a way that they also satisfy $\epsilon_\beta=1$ for every $\beta\in\Phi$.  This is achieved by replacing some of the $x_{\pm\gamma}(\xi)$ by $x'_{\pm\gamma}(\xi):=x_{\pm\gamma}(-\xi)$.  We shall choose such a parametrisation.    

For a subset $\Pi\subset \Delta$ we shall denote by $\Phi_\Pi$ the root system generated by $\Pi$ and 
by $P_\Pi$ the standard parabolic subgroup containing $B$ associated with $\Pi$, i.e., such that its standard Levi subgroup $L_\Pi$ is generated by $T$ and by the root subgroups corresponding to roots in $\Phi_\Pi$.  The intersection $U\cap L_\Pi$ will be denoted by $U_\Pi$.
If $\alpha\in\Delta$  then we shall put $P_\alpha$ to indicate $P_{\{\alpha\}}$ . For any parabolic subgroup $P$ of $G$ we will denote by $P^u$ its unipotent radical. 
The parabolic subgroup of $W$ generated by the simple reflections with respect to roots in $\Pi\subset \Delta$ will be denoted by $W_\Pi$. 

For a subgroup $H$ of $G$ we shall denote by $Z(H)$ its center and by  $H^\circ$ its identity component. When an automorphism $\tau$ acts on an algebraic structure $S$ (e.g. a subgroup of $G$ or a root system) we shall indicate by $S^\tau$ the substructure of $\tau$-invariant elements of $S$. 

\section{Twisted conjugacy classes}\label{twisted-classes}

A $\theta$-twisted conjugacy class in $G$ is an orbit for the $G$-action on itself by $g\cdot_\theta x=gx\theta(g)^{-1}$. When there is no ambiguity on $\theta$, we shall call it also a twisted conjugacy class and we shall use the simplified notation $g*x$ for $g\cdot_\theta x$. The $\theta$-stabilizer of $x\in G$ in a subgroup $H$ of $G$ is the stabilizer for the $*$-action and it will be denoted by $H_x$. 

Let $C$ be a twisted conjugacy class of $G$. Since $C$ is an irreducible variety there exists a unique element  in $W$ for which $C\cap BwB$ is dense
in $C$. We shall denote this element by $w_C$. We have
$$C\subset \overline{C}=\overline{C\cap B w_C B}\subset
\overline{B w_CB}=\bigcup_{\sigma\leq w_C}B\sigma B$$
so the element $w_C$ is the maximum among those $w\in W$ for which $BwB\cap C$ is non-empty
 (cfr. \cite[Section 1]{ccc}).
The collection of $B$-orbits for the $*$-action in $C$ will be denoted by $\V$.
We will call {\em maximal orbits} the elements $v$ in $\V$ lying in
$B w_CB$ and we shall denote by $\V_{\max}$ the set of maximal
$B$-orbits in $C$.

\begin{definition}An element $w\in W$ is called a {\em $\theta$-twisted involution} if $w\theta(w)=1$.  \end{definition}

If there is no ambiguity on the automorphism, we shall also say that $w$ is a twisted involution.
It is shown in \cite{lu} that $w_C$ is always a twisted involution and a genuine involution in $W$, that it commutes with the automorphism $\theta$ of $\Phi$ and with the longest element $w_0$ in $W$, and that  it is of the form $w_0w_\Pi$ where $\Pi$ is a suitable $\theta$-invariant subset of $\Delta$ and $w_\Pi$ is the longest element in $W_\Pi$. Although the general assumption in the paper is that the base field is of  characteristic zero, the arguments used in Section 3 from Lemma 3.1 until Proposition 3.7  are characteristic-free. The set $\Pi$ is recovered from $w_C$ by the equality
\begin{equation}\label{Pi}\Pi=\{\alpha\in\Delta~|~w_C\theta\alpha=\alpha\}.\end{equation}
The list of possible $w_C$'s for $\theta$ a non-trivial automorphism of the Dynkin diagram is provided in \cite[Proposition 3.7]{lu}. We report here the list of possible pairs $(\Phi, \Pi)$ for every non-trivial $\theta$ for completeness.
\begin{equation}\label{list}
\begin{array}{cl}
(\Phi,\emptyset) &\mbox{for any $\Phi$ and any $\theta$;}\\
(A_{2n+1}, \{\alpha_1,\,\alpha_3,\,\ldots,\,\alpha_{2n+1}\}) &\theta=-w_0;\\
(D_4, \{\alpha_2\}) &\theta^3=1;\\
(D_4, \{\alpha_2, \alpha_i,\theta\alpha_i\}) & \mbox{$\theta^2=1$ and $\alpha_i\neq\theta\alpha_i$;}\\
(D_{2n} ,\{\alpha_{2l},\alpha_{2l+1},\ldots, \alpha_{2n-1},\alpha_{2n}\})& \mbox{$n>2$,  $1\leq l\leq n-1$ and } \theta\alpha_{2n-1}=\alpha_{2n};\\
(D_{2n+1} ,\{\alpha_{2l},\,\alpha_{2l+1},\,\ldots \alpha_{2n},\alpha_{2n+1}\}) & \mbox{$n\geq 2$,  $1\leq l\leq n$, and  $\theta=-w_0$;}\\
(E_6,\{\alpha_2,\,\alpha_3,\,\alpha_4,\,\alpha_5\})&\theta=-w_0. \\
\end{array}
\end{equation}
By definition, maximal $B$-orbits in $C$ are contained in $Bw_CB$. By $B*$conjugacy, we can make sure that every such orbit contains an element of the form $\m\v$ where $\v\in U$ and $\m$ is a representative of $w_C$ in $N(T)$.  We analyze now the possible representatives of a $\theta$-twisted conjugacy class lying in a maximal $B$-orbit.
\begin{lemma}\label{minimal}
Let $C$ be a $\theta$-twisted conjugacy class and let $w_C=w_0w_\Pi$. Let $x=\m\v\in C\cap T w_C U$ for some lift $\m$ of $w_C$ in $N(T)$. Then $\v\in P_\alpha^u$ for every $\alpha\in\Pi$.
\end{lemma}
\pf Let $\v=x_\alpha(\xi)\v'$ for some $\v'\in P^u_\alpha$ and $\xi\in k$. Let $\s_\alpha\in N(T)$ be as in Section \ref{notation}.
We consider 
$$y=\theta^{-1}(\s_\alpha)*x=\theta^{-1}(\s_\alpha)\m x_{\alpha}(\xi)\v' \s_{\alpha}^{-1}=t\m \s_\alpha x_{\alpha}(\xi)\s^{-1}_\alpha \v''\in T\m x_{-\alpha}(\eta\xi)B$$ for some $t\in T$, $\eta\in k^*$ and $\v''\in P_\alpha^u$. Here we have used \eqref{Pi}.

If $\xi\neq0$, then  $y\in Bw_C Bs_\alpha B$. Since $w_C\alpha=\theta\alpha$ is a positive root, then $y\in C \cap B w_C s_\alpha B$ with $w_Cs_\alpha> w_C$ in the Bruhat order, a contradiction.\hfill$\Box$

\begin{lemma}\label{radical}
Let $C$ be a $\theta$-twisted conjugacy class and let $w_C=w_0w_\Pi$. Let $x=\m\v\in C\cap T w_C U$ for some  lift $\m$ of $w_C$. Then $\v\in P_\Pi^u$.
\end{lemma}
\pf We will show by induction on the height of $\beta\in\Phi_\Pi$ that once we fix an ordering of $\Phi^+$, the coefficient $c_\beta$ of $x_\beta$ in the expression of $\v$ as a product of elements in the root subgroups is trivial. We assume that the fixed ordering is compatible with the height of the roots. The basis of the induction is Lemma \ref{minimal}.  Let $\beta$ be the first root in $\Phi_\Pi$ for which $c_\beta\neq0$ and let its height be $h$. Then, $\v=\v_1 x_\beta(c_\beta)\v_2$ for some $\v_1\in P_\Pi^u$ and some $\v_2$ in a product of root subgroups  associated with roots of height greater or equal than $h$ and different from $\beta$. Since $\Phi$ is simply-laced, there exists $w\in W_\Pi$ such that $\ell(w)=h-1$ and $w\beta=\alpha\in\Pi$. Let $\w$ be a lift of $w$ in $N(T)$.
We consider the following representative of $C$:
$$\begin{array}{rl}
y&=\theta^{-1}(\w)*x=\theta^{-1}(\w)x\w^{-1}=\theta^{-1}(\w)\m \v_1 x_\beta(c_\beta)\v_2\w^{-1}\\
&=t\m(\w \v_1\w^{-1})(\w x_\beta(c_\beta)\w^{-1})(\w \v_2\w^{-1})=t\m \v'_1 x_\alpha(\eta c_\beta)\v'_2
\end{array}
$$
for some $t\in T$ and some $\eta\in k^*$. The element $\w$ normalizes $P_\Pi^u$ and by the hypothesis on the height also $\v_2'\in U$, so $\v_1',\,\v_2'\in P_\alpha^u$. Hence, $y\in T\m U$. By Proposition \ref{minimal} we necessarily have $c_\beta=0$.\hfill$\Box$

\begin{lemma}\label{orbit}Let $C$ be a $\theta$-twisted conjugacy class and let $\alpha\in \Phi^+$. Assume that  for every $x=\m\v\in C\cap T w_C U$ the coefficient of $x_\alpha$ in the expression of $\v$ as a product of elements in the root subgroups is zero for any ordering of the positive roots. Then, for every such $x$ the  coefficient of $x_\beta$ in the expression of $\v$ is zero for every  $\beta\in W_\Pi\alpha$ and for every ordering of the positive roots.
\end{lemma}
\pf If $\alpha\in\Phi_\Pi$ this is clear by Lemma \ref{radical} so we may assume $\alpha\in\Phi^+\setminus\Phi_\Pi$.

Let $\alpha=w\beta$ with $w\in W_\Pi$ and let $\w$ be a lift of $w$ in $N(T)$. We will write $\v=\v_1 x_{\beta}(c_\beta)\v_2$ for some $\v_1,\,\v_2\in P_\Pi^u$, products in root subgroups different from $X_\beta$. We consider the element
$$y=\theta^{-1}(\w)*x=\theta^{-1}(\w) x\w^{-1}=\theta^{-1}(\w)\m\w^{-1}\v'_1 x_{\alpha}(c_\beta')\v_2'$$
for some $c_\beta'\in k$ which is nonzero if and only if $c_\beta$ is nonzero. Since $W_\Pi\delta\in\Phi^+$ for every $\delta\in\Phi^+\setminus\Phi_\Pi$, we have $\v_1',\,\v_2'\in P_\Pi^u$. Moreover,  $w_C\theta\gamma=\gamma$ for every $\gamma\in \Phi_\Pi$ so $w_C^{-1}\theta^{-1} w \theta w_C=w$. Thus, $y\in C\cap T w_CU$ and by the assumption $c_\beta'=0$, whence the statement.
\hfill$\Box$

\begin{lemma}\label{levi-commute}Let $C$ be a $\theta$-twisted conjugacy class and let $x=\m\v\in Tw_CU\cap C$. Then,  
$[L_\Pi,L_\Pi]$ lies in the $\theta$-stabilizer of $\m$. \end{lemma}
\pf Let $\alpha\in\Pi$ and let $\beta=\theta\alpha\in \Pi$. We have $\theta x_{\alpha}(\xi)=x_\beta(\xi)$.  By Lemma \ref{minimal} we know that $\v\in P_\alpha^u$. We consider the following representative of $C$:
$$\begin{array}{rl}
y&=x_\alpha(\xi)x \theta(x_\alpha(-\xi))=x_\alpha(\xi)\m\v x_\beta(-\xi)=\m x_\beta(\eta \xi)\v x_{\beta}(-\xi)
\end{array}
$$
for some $\eta\in k^*$.  Here we have used that $w_C\theta\alpha=\alpha$. By Lemma \ref{minimal}, we have $x_\beta(\eta \xi)\v x_{\beta}(-\xi)\in P_\beta^u$ 
and  this is possible only if $\eta=1$, that is, if the root subgroup $X_\alpha$ lies in the $\theta$-stabilizer of $\m$.

Let us now consider $-\alpha$ and $-\beta$. Again we have $\theta x_{-\alpha}(\xi)=x_{-\beta}(\xi)$. We consider the following representative of $C$:
$$\begin{array}{rl}
z&=x_{-\alpha}(\xi)x x_{-\beta}(-\xi)=x_{-\alpha}(\xi)\m\v x_{-\beta}(-\xi)\\
&=\m x_{-\beta}(\eta \xi)\v x_{-\beta}(-\xi)\in\m x_{-\beta}(\eta \xi-\xi)P_\beta^u
\end{array}
$$
for some $\eta\in k^*$. If we had $\eta\neq1$ we would have
$z\in \m B s_\beta B\subset B w_Cs_\beta B$ because $w_C\beta=\theta^{-1}\beta\in\Phi^+$. This would contradict maximality of $w_C$, hence $\eta=1$ and the root subgroup $X_{-\alpha}$ lies in the $\theta$-stabilizer of $\m$. 
 \hfill$\Box$

\bigskip

When $\theta^2=1$, the class $G*1=\{g\theta(g)^{-1}~|~g\in G\}$, extensively studied in   \cite{R,results, RS},  intersects only Bruhat cells corresponding to twisted involutions in $W$ because $\theta(g\theta(g)^{-1})=(g\theta(g)^{-1})^{-1}$. We are going  to study all $\theta$-twisted conjugacy classes sharing this property.

\begin{definition}Let $C$ be a $\theta$-twisted conjugacy class. We will say that $C$ is {\em involutive} if $C\cap BwB\neq\emptyset$ only when $w$ is a twisted involution.
\end{definition}

\begin{remark}\label{counter}{\rm If $\Phi$ is of type $D_4$ and  $\theta$ is the automorphism of order $3$ mapping $\alpha_1$ to $\alpha_3$, then there exist no involutive $\theta$-twisted conjugacy classes. Indeed, given a representative $x=\m\v\in C\cap Tw_CU$, if $w_C=w_0s_2$ then the representative $y=\s_1* x$ lies in $Bs_1w_0s_2s_3B\cup Bs_1w_0s_2B$ and both Weyl group elements are not twisted involutions. 
If instead $w_C=w_0$, then $\s_1* x\in Bs_1w_0s_3B\cup Bs_1w_0 B$ and we conclude as above.}
\end{remark}

\section{Spherical twisted conjugacy classes}\label{spherical-classes}

In this section we will introduce spherical $G$-spaces and we will show that if $\theta$ is an involution, then every spherical $\theta$-twisted  conjugacy class is involutive. 

\begin{definition}A transitive $G$-variety is called {\em spherical} if it has a dense $B$-orbit.
\end{definition}

The dense $B$-orbit is necessarily unique. It has been shown in \cite{Bri,Vin,knop}  that a $G$-variety is spherical if and only if $B$ acts on it with finitely many orbits. 

The following Lemma is a $\theta$-twisted analogue of \cite[Lemma 3.1]{gio-fourier}. We report the proof to keep the paper self contained.
 
 \begin{lemma}\label{finite}Let $C$ be a $\theta$-twisted conjugacy class. The following are equivalent
 \begin{enumerate}
 \item $C$ is spherical.
 \item $\V_{\rm max}$ is a finite set.
  \end{enumerate} 
 \end{lemma}
 \pf  One implication is immediate from the above remarks. We have
 $$\overline{C}=\overline{C\cap B w_C B}=\overline{\bigcup_{v\in{\cal V}_{\rm max}} v}.$$ If $\V_{\rm max}$ is finite, then irreducibility of $C$ forces $\overline{v}=\overline{C}$ for some $v\in\V_{\rm max}$.
 \hfill$\Box$
 
 \bigskip
 
Let $M(W)$ denote the monoid generated by the symbols $r_\alpha$ for $\alpha\in\Delta$ subject to the braid relations and the relation $r_\alpha^2=r_\alpha$. 
Given a spherical $G$-variety,  there are an $M(W)$-action and a $W$-action on the set of its $B$-orbits  ${\cal V}$. These actions have been introduced in \cite{RS} and \cite{knop}, respectively, and they have been further analyzed and applied in  \cite{closures}, \cite[\S 4.1]{MS}, \cite{springer-schubert}. 
For  $v\in\V$, the $B$-orbit $r_\alpha(v)$ is the dense $B$-orbit in $P_\alpha v$.  In order to introduce the $W$-action we need to provide more background information.

Let $v\in \V$. Then, the action of $P_\alpha$  on $P_\alpha/B\cong {\mathbb P}^1$ defines a group morphism $\psi\colon P_\alpha\to PGL_2(k)$ whose kernel is ${\rm  Ker}(\alpha)P^u_\alpha$. The stabilizer $(P_\alpha)_y$  of a point $y$ in $v$ acts on $P_\alpha/B$ with finitely
many orbits.  
The image $H$ of $(P_\alpha)_y$ in $PGL_2(k)$ is of one of the following types:
$PGL_2(k)$; solvable and contains a connected nontrivial unipotent subgroup;
a torus;  the normalizer of a torus. More precisely, we fall in one of the following cases:
\begin{enumerate}
\item[I] $P_\alpha v=v$ and $H=PGL_2(k)$; 
\item[II]$P_\alpha v=v\cup v'=P_\alpha v'$ where  $v'= r_\alpha(v)$ or $v=r_\alpha(v')$, with $|\dim v-\dim v'|=1$ and  $H$ is solvable containing a nontrivial unipotent subgroup.
 \item[III] $P_\alpha v=v\cup v'\cup v'' = P_\alpha v'=P_\alpha v''$, with $v,v',v''$ distinct, where either $v=r_\alpha(v')=r_\alpha(v'')$ and $\dim v=\dim v'+1=\dim v''+1$, or $v'=r_\alpha(v)=r_\alpha(v'')$ and $\dim v'=\dim v+1=\dim v''+1$,   and $H$ is a torus.
 \item[IV]  $P_\alpha v=v\cup v'=P_\alpha v'$ where $v'= r_\alpha(v)$ or $v=r_\alpha(v')$, with $|\dim v-\dim v'|=1$ and $H$ is the normalizer of a torus.
\end{enumerate}
The $W$-action on $\V$ can be defined as follows  (\cite{knop}, \cite[\S 4.2.5,
  Remark]{MS}):  the simple reflection $s_\alpha$ interchanges 
  the two $B$-orbits in case II; it interchanges the  two non-dense orbits in case III and it fixes all $B$-orbits in types I and IV and the dense $B$-orbit in type III. 
The action of $s_\alpha$ on $v$ will be denoted by $s_\alpha.v$. 

We will now show that every $v\in{\cal V}$ can be reached from a closed one by means of a path in which each step is given either by the action of $s_\alpha\in W$ or the action of $r_\beta\in M(W)$.  This is formalized as follows.
\begin{definition}(\cite[\S 3.6]{springer-schubert}) Let $X$ be a spherical  $G$-variety and let $\V$ be its set of $B$-orbits. A {\em reduced decomposition}
of $v\in \V$ is a pair $({\bf v},{\bf s})$ with ${\bf v}=(v(0),v(1),\ldots, v(r))$ a sequence of
 elements in $\V$ and ${\bf s}=(s_{i_1},\ldots, s_{i_r})$ a
sequence of simple reflections such that: $v(0)$ is closed;
$v(j)=r_{i_j}(v({j-1}))$ for $1\leq j\leq r$;
$\dim(v(j))=\dim(v(j-1))+1$ and $v(r)=v$.
\end{definition}

\begin{remark}\label{inizio}{\rm Every $B$-orbit  $v$ in a symmetric space admits a reduced
decomposition  by \cite[\S 7]{RS}.  However, it is shown in \cite[Example 2]{closures} that this is not always the case for a spherical $G$-variety $X$.
On the other hand, for every closed $B$-orbit $v$ in $X$, there exists a reduced decomposition of the dense $B$-orbit $v_0$ with $v(0)=v$. 
Indeed, if $r_\alpha v'\neq v'$ for $v'\in\V$,  then $\dim r_\alpha v'=\dim v'+1$,  so we may inductively construct a sequence $(v(0),v(1),\ldots, v(r))$ with $v(0)=v$ and $v(j)=r_{i_j}(v({j-1}))$ for $1\leq j\leq r$.  We can choose  $s_{i_j}$ so that it satisfies $\dim(v(j))=\dim(v(j-1))+1$ provided there is $\alpha\in\Delta$ such that  $r_{\alpha}v(j-1)\neq v(j-1)$. The procedure will stop at some $B$-orbit $v(r)$ such that $r_\alpha(v(r))=v(r)$ for every $\alpha\in\Delta$. Then,
$$\overline{v(r)}=\overline{P_\alpha v(r)}=P_\alpha \overline{v(r)}$$
where we have adapted the argument in \cite[Exercise 6.2.11(5)]{springer}. Thus, $X=G\cdot\overline{v(r)}=\overline{v(r)}$  and therefore
 $v(r)=v_0$. The same argument shows that any sequence $(v(0),v(1),\ldots, v(r))$ with $v(j)=r_{i_j}(v({j-1}))$ for $1\leq j\leq r$ and 
$\dim(v(j))=\dim(v(j-1))+1$ can be completed to a reduced decomposition of the unique dense $B$-orbit. }
\end{remark}

A weaker notion of reduced decomposition exists for every $v\in \V$. 

\begin{definition}\label{sub}(\cite[\S 3.6]{springer-schubert}) A  {\em subexpression} of 
a reduced decomposition $({\bf v},{\bf s})=((v(0),\ldots, v(r)),(s_{i_1},\ldots, s_{i_r}))$ of $v\in \V$  is a sequence 
${\bf  x}=(v'(0),v'(1),\ldots,\,v'(r))$ of elements in $\V$ with
 $v'(0)=v(0)$ and such that for $1\leq
j\leq r$ only one of the following alternatives occurs:
\begin{enumerate}
\item[(a)] $v'(j-1)=v'(j)$; 
\item[(b)] $\dim v'(j-1)=\dim v'(j)-1$ and $v'(j)=r_{i_j}(v'(j-1))$;
\item[(c)] $v'(j-1)\not=v'(j)$, $\dim v'(j-1)=\dim v'(j)$ and $v'(j)=s_{i_j}.v'(j-1)$.
\end{enumerate}
The element $v'(r)$ is called the {\em final term} of the subexpression.
\end{definition}
By \cite[\S 3.6 Proposition
  2]{springer-schubert},  if $v'\in \V$ has a reduced decomposition $({\bf v'},{\bf s})$, then for any $v\in\V$ which is contained in the closure of $v'$, there exists a subexpression of $({\bf v'},{\bf s})$ with end-point $v$. In particular, for every $v\in \V$,  there exists a subexpression of any reduced decomposition of the
dense $B$-orbit $v_0$ admitting $v$ as final term. The statement is given in characteristic zero but the proof holds also in positive odd characteristic.

\smallskip

The following theorem is a generalization of \cite[Theorem 2.7]{gio-mathZ}, from the untwisted case. The argument has been shortened and simplified and it  also works for $\theta$  trivial.

\begin{theorem}\label{induction}Let $\theta$ be an involution of the Dynkin diagram of $G$ and let $C$ be a spherical $\theta$-twisted conjugacy class.
Then $C$ is involutive.
\end{theorem}
\pf By \cite[Lemma 7.3]{steinberg-endomorphisms} we may choose a representative $y\in B$ for $C$. Hence, $B* y\subset B$ and its closure contains a closed $B$-orbit $x(0)$ lying in $B$. By Remark \ref{inizio} there is a reduced decomposition $({\bf v_0}, {\bf s})$ of the dense $B$-orbit  $v_0$ with initial point $x(0)$. Let $v\in\V$. By \cite[\S 3.6 Proposition
  2]{springer-schubert}, there is a subexpression ${\bf  x}=(v'(0),v'(1),\ldots,\,v'(r))$
 of $({\bf v_0}, {\bf s})$ with initial point $v'(0)=x(0)$ and final point $v'(r)=v$.

We will show by induction on $j$ that $v'(j)$ lies in the Bruhat cell corresponding to a twisted involution. For $j=0$ this is immediate. Let us assume that $v'(j-1)\subset Bw_{j-1}B$ with $w_{j-1}$ a twisted involution. We consider the step from $v'(j-1)$ to $v'(j)$. 
If we are in case (a) of Definition \ref{sub} there is nothing to prove. 
If we are in case (b) let $\alpha=\alpha_{i_j}$. Then $P_\alpha* v'(j-1)\subset Bw_{j-1} B\cup Bs_\alpha Bw_{j-1}Bs_{\theta\alpha}B.$
According to \cite[Lemma 3.2]{results} there are three possibilities: 
\begin{itemize}
\item $\ell(s_\alpha w_{j-1}s_{\theta\alpha})=\ell(w_{j-1})+2$ so 
$P_\alpha *v'(j-1)\subset Bw_{j-1} B\cup Bs_\alpha w_{j-1}s_{\theta\alpha}B$ and both Weyl group elements involved are twisted involutions;
\item $s_\alpha w_{j-1}=w_{j-1}s_{\theta\alpha}$ so
$P_\alpha *v'(j-1)\subset Bw_{j-1} B\cup Bs_\alpha w_{j-1}B$ and both Weyl group elements involved are twisted involutions;
\item $\ell(s_\alpha w_{j-1}s_{\theta\alpha})=\ell(w_{j-1})-2$ so
$$P_\alpha* v'(j-1)\subset Bw_{j-1} B\cup Bs_\alpha w_{j-1}B\cup Bw_{j-1}s_{\theta\alpha}B\cup Bs_\alpha w_{j-1}s_{\theta\alpha}B.$$ 
Since $v'(j)=r_{\alpha}(v'(j-1))$ is dense in $P_\alpha* v'(j-1)$, it lies in a cell corresponding to a $\sigma\in W$ with $\sigma\geq w_{j-1}$ in the Bruhat order. Thus, 
$v'(j)\subset Bw_{j-1} B$.
\end{itemize}
If we are in case (c) then we are necessarily in situation III in the description of the $W$-action on $\V$ and $v'(j-1),\,v'(j)$ are the non-dense $B$-orbits in 
$P_\alpha *v'(j-1)=P_\alpha *v'(j)$. If  $\ell(s_\alpha w_{j-1}s_{\theta(\alpha)})=\ell(w_{j-1})+2$ or if $s_\alpha w_{j-1}=w_{j-1}s_{\theta\alpha}$ we may proceed as in case (b). 
Let us assume that $\ell(s_\alpha w_{j-1}s_{\theta\alpha})=\ell(w_{j-1})-2$ so
$$P_\alpha* v'(j-1)\subset Bw_{j-1} B\cup Bs_\alpha w_{j-1}B\cup Bw_{j-1}s_{\theta\alpha}B\cup Bs_\alpha w_{j-1}s_{\theta\alpha}B.$$  
Let $x_1\in Tw_{j-1}U\cap v'(j-1)$, and $x_2\in U w_{j-1} T\cap v'(j-1)$. We have 
$$y_1:=\s_\alpha* x_1\in P_\alpha *v'(j-1)\cap (Bs_\alpha w_{j-1}B\cup Bs_\alpha w_{j-1}s_{\theta\alpha}B)$$
$$y_2:=\s_\alpha* x_2\in P_\alpha* v'(j-1)\cap (Bs_\alpha w_{j-1}s_{\theta\alpha}B\cup B w_{j-1}s_{\theta\alpha}B).$$
Thus $y_1,\,y_2\in v'(j)$ because there are only three $B$-orbits in $P_\alpha *v'(j-1)$ and  by the discussion of case (b) we have $r_\alpha(v'(j-1))\subset Bw_{j-1}B$.  Hence, $v'(j)\subset Bs_\alpha w_{j-1}s_{\theta\alpha}B$ and $s_\alpha w_{j-1}s_{\theta\alpha}$ is a twisted involution.\hfill$\Box$

\smallskip

\begin{remark}{\rm Theorem \ref{induction} fails if we drop the assumption on $\theta$ to be an involution. Indeed, in the triality case  it has been shown in \cite[Example 3.9]{lu} and \cite[Section 4.5]{mauro}  that the class $G*1$ is spherical. However,  it is not involutive by Remark \ref{counter}.}
\end{remark}

\section{Involutive twisted conjugacy classes}\label{involutive}

This section is devoted to the understanding of involutive $\theta$-twisted conjugacy classes so we shall assume that $\theta$ is an involution. We aim at getting some control on the representatives of $C$ in maximal $B$-orbits. 

\begin{lemma}\label{simple}Let $C$ be an involutive $\theta$-twisted conjugacy class, let $w_C=w_0w_\Pi$ with $\Pi\neq\emptyset$ and let $x=\m\v\in C\cap T w_C U$. Then $\v\in P_\alpha^u$ for every $\alpha\in\Delta$ such that $\alpha\not\perp \Pi$.
\end{lemma}
\pf If $\alpha\in\Pi$ this is Lemma \ref{minimal}. Let $\alpha\in\Delta\setminus \Pi$, let $\v=x_{\alpha}(c)\v'$ with $\v'\in P_\alpha^u$. We consider
$$\begin{array}{rl}
y&=\theta^{-1}(\s_\alpha)*x=\theta^{-1}(\s_\alpha)x \s_\alpha^{-1}=\theta^{-1}(\s_\alpha)\m  \s_\alpha^{-1} x_{-\alpha}(c')\v''
\end{array}
$$ for some $\v''\in P_\alpha^u$ and some $c'\in k$ which is nonzero if and only if $c$ is nonzero. If $c'\neq0$ then
$y$ lies in $T s_{\theta^{-1}\alpha} w_C s_\alpha B s_\alpha B\cap C$. It follows from a straightforward verification that if $\alpha\not\perp\Pi$ we have $w_\Pi\alpha\in\Phi^+\setminus\Delta$, so $\beta=w_C\alpha\in-(\Phi^+\setminus\Delta)$. Thus, $s_{\theta^{-1}\alpha}w_C s_\alpha\alpha\in\Phi^+$ and $y\in C\cap Bs_{\theta^{-1}\alpha}w_CB$. However, $s_{\theta^{-1}\alpha}w_C$ is not a twisted involution. Indeed,
$$s_{\theta^{-1}\alpha}w_C\theta(s_{\theta^{-1}\alpha}w_C)=s_{\theta^{-1}\alpha}(w_C s_\alpha w_C)=s_{\theta^{-1}\alpha}s_\beta\neq1$$ where we have used that $w_C$ is $\theta$-invariant. Hence, $c'=c=0$ and we have the statement. \hfill$\Box$

\bigskip

In the spirit of \cite{results} we define the following subsets of roots for $w=w_C=w_0w_\Pi\in W$ a Weyl group element in the list \eqref{list}.
$$\begin{array}{rl}
C_w&=\{\alpha\in\Phi^+~|~w\theta\alpha\in -\Phi^+\mbox{ and }w\theta\alpha\neq-\alpha\},\\
I_w&=\{\alpha\in\Phi^+~|~w\theta\alpha=\alpha\},\\
R_w&=\{\alpha\in\Phi^+~|~w\theta\alpha=-\alpha\}.
\end{array}
$$
Such sets are called the set of complex, imaginary and real roots relative to $w$, respectively.

Since for $\alpha\in\Phi^+$ we have  $w_\Pi\alpha\in-\Phi^+$ if and only if $\alpha\in\Phi_\Pi$, we have $I_w=\Phi_\Pi\cap\Phi^+$. 
Besides, $\Phi^+$ is the disjoint union of $I_w$, $R_w$, and $C_w$.

The set $R_w$ is contained in the $(-1)$-eigenspace of the orthogonal map $w\theta$, therefore it lies in $I_w^\perp=\Pi^\perp$ so, for every $\alpha\in R_w$ we have $w_0\theta\alpha=-\alpha$.
On the other hand, if $\beta\in\Pi^\perp\cap\Phi^+$ and $w_0\theta\beta=-\beta$
then $w\theta\beta=\theta w_0w_\Pi\beta=\theta w_0\beta=-\beta$. Hence, we have
$$
R_w=\left\{\begin{array}{ll}
\Pi^\perp\cap\Phi^+&\mbox{if $w_0=-\theta$,}\\
(\Pi^\perp\cap\Phi^+)^\theta&\mbox{if $w_0=-1$.}\
\end{array}\right.
$$

The union $R=R_{w}\cup (-R_w)$ is a root subsystem of $\Phi$ and we may consider the reductive subgroup $G_{R}=\langle T,\,X_\alpha~|~\alpha\in R\rangle$. We may choose a set $\Delta_R=\{\gamma_1,\,\ldots,\,\gamma_r\}\subset \Phi^+$ of simple roots in $R_{w}$ so that $B\cap G_R$ is a Borel subgroup of $G_R$ and $U_R=U\cap G_R$ is its unipotent radical. 
The subset $\Delta_R$ need not be a subset of $\Delta$.

The Weyl group $W_R$ of $G_R$ is generated by some reflections in $W$ so it is a subgroup of $W$. 

\begin{remark}\label{slodowy}{\rm Since the root system $R$ of $G_R$ is ${\mathbb Q}$-closed, it follows from \cite[Section 3.5]{slodowy} that $G_R$ is the Levi factor of a parabolic subgroup of $G$. Hence, its derived subgroup is simply-connected and $W_R$ is conjugate to a parabolic subgroup of $W$. }\end{remark}

\begin{proposition}\label{real}With the above notation, let $C$ be an involutive $\theta$-twisted conjugacy class in $G$ and let $x=\m\v\in T w_CU$. Then $\v$ lies in $U_R$.
\end{proposition}
\pf If $\theta=-w_0$ and $w_C=w_0$ (i.e. $\Pi=\emptyset$) this condition is empty.  

The basic idea of the proof for all non-trivial cases is to exhibit, for $\beta\in C_{w_C}$,  a Weyl group element $\sigma$ satisfying the following properties:
\begin{enumerate}
\item[(1)] $\alpha=\sigma\beta\in\Delta$;
\item[(2)] $\sigmad \v \sigmad^{-1}\in U$ for a representative $\sigmad\in N(T)$;
\item[(3)] the root $\gamma=\theta^{-1}(\sigma)w_C\beta$ lies in $-(\Phi^+\setminus\{\theta^{-1}\alpha\}).$
\end{enumerate}

Then, the element
$y=\theta^{-1}(\sigmad)*x$ lies in $T\theta^{-1}(\sigma)w_C\sigma^{-1} x_\alpha(\eta c_\beta)P_\alpha^u$ for some $\eta\in k^*$. Hence, if the coefficient $c_\beta$ of $x_\beta$ in the expression of $\v$ is non-zero we have 
$$z=\theta^{-1}( \s_\alpha)*y\in  Ts_{\theta^{-1}\alpha}\theta^{-1}(\sigma)w_C\sigma^{-1}s_\alpha^{-1}B s_\alpha B.$$
As $\gamma$ lies in $-\Phi^+$ and it is different from $-\theta^{-1}\alpha$, we have, for  $\tau=s_{\theta^{-1}\alpha}\theta^{-1}(\sigma)w_C\sigma^{-1}$ the inequality $\tau>\tau s_\alpha$ so
$z$ lies in $B\tau B\cap C$. However, 
$$\begin{array}{rl}
\tau\theta(\tau)&=s_{\theta^{-1}\alpha}\theta^{-1}(\sigma)w_C\sigma^{-1}s_\alpha\sigma w_C\theta (\sigma^{-1})\\
&=s_{\theta^{-1}\alpha}(\theta^{-1}(\sigma)w_C s_\beta w_C^{-1} \theta (\sigma^{-1}))\\
&=s_{\theta^{-1}\alpha}s_\gamma\neq 1.
\end{array}$$
Therefore, if $c_\beta$ is non-zero then $C$ is not involutive.

We discuss the different cases separately, according to the classification of the $w_C$'s in \eqref{list}.

\smallskip

\noindent{\bf Case $(A_{2n+1}, \{\alpha_1,\alpha_3,\,\ldots,\,\alpha_{2n+1}\})$.} We will show that $C$ is the twisted conjugacy class of a lift of $w_C$. Let $x=\m\v\in T w_CU\cap C$ and
let us assume that $\v=\prod x_{\gamma}(c_\gamma)$ with $c_\gamma=0$ for $\gamma$ of height smaller than $h$.  Then $h\geq 2$ by Lemma \ref{minimal} and Lemma \ref{simple}. Let  $\beta=\alpha_i+\cdots+\alpha_{j}$ be a root of minimal height for which $c_\beta\neq0$.
If $j$ or $i$ were odd then we could apply Lemma \ref{orbit} obtaining a contradiction, so $i$ and $j$ are even. Then
the Weyl group element $\sigma=s_{j-1}s_{j-2}\cdots s_i$ satisfies properties (1), (2), (3) and we have the statement in this case.


\noindent{\bf Case $(D_n, \{\alpha_j\}_{j\geq 2l})$, for $n\geq4$ and $1\leq l\leq\left[\frac{n-1}{2}\right]$.} Let $\theta$ be the automorphism of the Dynkin diagram of type $D_n$ interchanging $\alpha_n$ and $\alpha_{n-1}$. If $n$ is even then $w_0=-1$ whereas if $n$ is odd $w_0=-\theta$. 
If $C$ is an involutive $\theta$-twisted conjugacy class with $w_C=w_0w_\Pi$, for $\Pi=\{\alpha_j\}_{j\geq 2l}$ with $1\leq l\leq m-1$ if $n=2m$ and $1\leq l\leq m$ if $n=2m+1$, we shall show that
the coefficient of $x_\beta$ in the expression of $\v$ is trivial for every $\beta$ which is not orthogonal to $\Pi$.
If $\beta$ is simple this is Lemma \ref{simple}. 
By Lemma \ref{orbit} it is enough to prove the statement for the roots of the form $\beta_{j}=\alpha_j+\cdots+\alpha_{2l-1}$.
Let $\beta=\beta_i$ be the root of minimum height among the $\beta_j$'s for which the coefficient is non-zero. Then $\beta$ is not simple and $\sigma=s_{i+1}\cdots s_{2l-1}$ satisfies properties (1), (2), (3), so $c_{\beta_j}=0$ for every $j$. The remaining cases in type $D_4$ follow by symmetry. 

\noindent{\bf Case $(D_{2n},\emptyset)$.} In this case the positive roots that are not $\theta$-invariant are of the form $\beta_i=\alpha_i+\cdots+\alpha_{2n-2}+\alpha_{2n-1}$ and $\theta\beta_i$ for $i\leq 2n-1$. 
Let $\beta=\beta_j$ be the root of minimal height of this form for which the coefficient in the expression of $\v$ is non-zero. Then $\beta$ is not simple and $\sigma= s_{j+1}\cdots s_{2n-2}s_{2n-1}$ satisfies properties (1), (2) and (3). Hence,  the coefficient of $\beta_i$ is zero for every $i$.
The case of $\theta\beta_j$ is treated similarly.

\noindent{\bf Case $(E_6, \{\alpha_2,\alpha_3,\alpha_4,\alpha_5\})$.} In this case $C$  is represented by a lift of $w_C$ in $N(T)$. Indeed, it follows from Lemma \ref{minimal}, Lemma \ref{simple} and Lemma \ref{orbit} that $\v$ can be expressed as a product in the root subgroups associated with the positive roots outside $\Phi_\Pi$, $W_\Pi\alpha_1$ and $W_\Pi\alpha_6$, that is, the positive roots in the orbit $W_\Pi\beta$ for $\beta=\alpha_1+\alpha_3+\alpha_4+\alpha_5+\alpha_6$. All positive roots in $W_{\Pi}\beta\setminus\{\beta\}$ have height strictly greater than $5$.
Then, the Weyl group element $\sigma=s_5s_4s_3s_1$ satisfies properties (1), (2) and (3) for the root $\beta$ so $\v=1$.
This concludes the proof of Proposition \ref{real}.\hfill$\Box$

\bigskip

An element in $W$ is called a twisted-identity (\cite{hult1}) if it is of the form $w\theta (w)^{-1}$ for some $w\in W$. A $\theta$-twisted conjugacy class is called $\theta$-semisimple if it has a representative in $T$ (\cite{springer-twisted}). 

\begin{corollary} Let $C$ be an involutive $\theta$-twisted conjugacy class such that  $w_C$ falls in one of the following cases: $(A_{2n+1}, \{\alpha_1,\alpha_3,\ldots,\,\alpha_{2n+1}\})$, $(D_n, \{\alpha_2,\,\ldots,\alpha_n\})$, $(E_6, \{\alpha_2,\alpha_3,\alpha_4,\alpha_5\})$. Then $C$ is $\theta$-semisimple.
\end{corollary}
\pf In these cases $C$ may be represented by some $\m\in w_C T$. It is enough to show that $w_C$ is a twisted identity because if $w_C=w\theta(w^{-1})$ for some $w\in W$, then we have $\w^{-1}* \m \in T\cap C$  for every lift $\w$ of $w$ in $N(T)$.

In type $A_{2n+1}$ the element $w_C=w_0w_\Pi$ is the permutation on $2n+2$ letters with cyclic decomposition 
$(1\  2n+1)(2\  2n+2)(3\  2n-1)(4\  2n)\cdots (n\  n+2)(n+1\  n+3)$ for $n$ odd and 
$(1 \ 2n+1)(2\  2n+2)(3\  2n-1)(4\  2n)\cdots (n-1\  n+3)(n\  n+4)$ with $n+1$ and $n+2$ fixed for $n$ even. In both expressions, each transposition $(a\  b)$ is followed by $\theta((a\  b))$ and since all transpositions in these expressions commute, $w_C=w\theta (w)^{-1}$ is always a twisted identity. 
In type $D_n$ one may verify that $w_C=w\theta (w)^{-1}$ for $w=s_{\alpha_1+\alpha_2+\cdots+\alpha_{n-1}}$.
In type $E_6$ we have $w_C=w\theta (w)^{-1}$ for $w=s_{\alpha_1+\alpha_2+\alpha_3+2\alpha_4+2\alpha_5+\alpha_6}$.\hfill$\Box$
 
\section{The intersections $C\cap \m U$}\label{quasi}

The aim of this section is to show that, unless we are in type $D_{2n}$ and $w_C=w_0$, every involutive $\theta$-twisted conjugacy class is spherical. The crucial step is Lemma \ref{sufficient}, where we conclude that it is enough to show that  that the number of suitable representatives of a maximal $B$-orbit in an involutive class is nonzero and finite. 

\begin{lemma}\label{sufficient} Let $C$ be a $\theta$-twisted conjugacy class and  let $\m\in N(T)$ be such that $C\cap \m U\neq\emptyset$. If 
\begin{enumerate}
\item for every $v\in \V_{\rm max}$ there is $z\in Z(G)$ such that $v\cap \m z U\neq\emptyset$,
\item $|C\cap \m z U|$ is finite for every $z\in Z(G)$
\end{enumerate}
then $C$ is spherical.
\end{lemma}
\pf Under the above assumptions we have:
$$
\begin{array}{rl}
|\V_{\rm max}|&=\sum_{v\in \V_{\rm max}}1\leq \sum_{v\in \V_{\rm max}}\sum_{z\in Z(G)}|v\cap \m zU|\\
&=\sum_{z\in Z(G)}|\bigcup_{v\in\V_{max}} v\cap \m z U|=\sum_{z\in Z(G)}|C\cap \m zU|<\infty 
\end{array}
$$
where we used that $Z(G)$ is finite. We conclude by using Lemma \ref{finite}.\hfill$\Box$

\bigskip

Let $C$ be an involutive $\theta$-twisted conjugacy class and let $\m\in N(T)$ be such that $C\cap \m U\neq\emptyset$. Then, for every $v\in\V_{\rm max}$ there is $x\in \m t U\cap v$ for some $t\in T$.
It follows from Lemma \ref{levi-commute} that  $[L_\Pi,L_\Pi]$ fixes $\m$ and $\m t$ under the $\theta$-action, and it is easy to conclude that it centralizes $t$.

\begin{proposition}\label{meets}Let $C$ be a $\theta$-twisted conjugacy class and let $\m U\cap C\neq \emptyset$. Let us assume that, for $w_C$, we are not in case $(D_{2m},\emptyset)$. Then, condition 1 in Lemma \ref{sufficient} is satisfied. 
\end{proposition}
\pf Let us consider the morphism 
$$\begin{array}{rl}
\psi\colon T&\to T\\
s&\mapsto (\m^{-1} s \m)\theta(s^{-1}).
\end{array}
$$
Then $\psi$ is a group morphism so its image is closed (\cite[Proposition 2.2.5]{springer}) and it lies in $Z(L_\Pi)$ by \eqref{Pi}.
It is also connected, so it lies in $Z(L_\Pi)^\circ$.
We recall that $\dim Z(L_\Pi)^\circ={\rm rk G}-|\Pi|$. 
On the other hand, ${\rm Ker}\psi=T^{w_C\theta}$.  By a simple direct computation, we see that, for all cases except from $(D_{2m},\emptyset)$, 
we have $\dim Z(L_\Pi)^\circ=\dim T-\dim T^{w_C\theta}$ so in all those cases ${\rm Im}\,\psi=Z(L_\Pi)^\circ$. 

Let $v\in \V_{\rm max}$ and let $x=\m t \v\in v$. Then for every $s\in T$ we have 
$$s* x=\m t\psi(s) \theta(s)\v\theta (s)^{-1}\in v\cap \m t\psi(s)U.$$ 
Thus, for every $r\in {\rm Im}(\psi)=Z(L_\Pi)^\circ$ we have  $v\cap \m t rU\neq\emptyset$. In the adjoint quotient of $G$ the center of any Levi factor of a parabolic subgroup is connected, so in $G$ we have $Z(L_\Pi)=Z(G)Z(L_\Pi)^\circ$ and $t$ lies in $z Z(L_\Pi)^\circ=z{\rm Im}(\psi)$ for some $z\in Z(G)$, whence the statement.\hfill$\Box$

In the following Lemmas we shall prove that if $C$ is involutive  then $|C\cap \m U|$ is finite for any $\m\in N(T)$.

\begin{lemma}\label{alpha-beta}Let $C$ be an involutive twisted conjugacy class. Let $\m$ be a representative of $w_C$ for which $\m U\cap C\neq\emptyset$.
Let $x=\m\v\in C\cap\m U$, with
  $\v=\prod_{\gamma\in R_{w_C}}x_\gamma(c_\gamma)$ in a fixed ordering of $R_{w_C}$. Let $\alpha$ and $\beta$ be adjacent simple
  roots in $\Delta_R$. Then, the number of possibilities for
  $c_\alpha$ and $c_\beta$ is finite. Moreover,  there is a $k_{\alpha,\beta}$ in $k$ depending only on the fixed ordering of the roots, on the structure constants of $G$, and on $\m$, such that $c_{\alpha+\beta}=k_{\alpha,\beta}c_\alpha c_\beta$.
\end{lemma}
\pf Let $P=P_{\{\alpha,\beta\}}$ be the standard parabolic subgroup of $G_R$  with unipotent radical $P^u$. Let us assume that $\alpha$ precedes $\beta$ in the
ordering of the roots in $R_{w_C}$. We may write: $x=\m\v\in\m
x_{\alpha}(c_\alpha)x_{\beta}(c_\beta)x_{\alpha+\beta}(c_{\alpha+\beta})P^u$. Let $\s_\gamma=x_\gamma(1)x_{-\gamma}(-1)x_\gamma(1)$ as in Section \ref{notation}, for $\gamma\in\{\alpha,\beta\}$.

The strategy of the proof is as follows: first we will show that for two precise values $h_1$ and $h_2$ of $h\in k$, depending on $c_\alpha$, the ordering, and the structure constants of $G$, the element  $y(h):=\theta^{-1}(\s_\alpha x_{\alpha}(h))*x$ lies in $B w_C s_\alpha B$. Then, we will consider the elements $\theta^{-1}(\s_\beta)*y(h_i)$ for $i=1,2$ and we will detect the Bruhat double cosets containing them. Imposing that this corresponds to a twisted involution will provide alternative necessary conditions on $c_{\alpha+\beta}$, $c_\beta$, $h_1$ and $h_2$. Then we will repeat the procedure interchanging the role of $\alpha$ and $\beta$, obtaining new alternative necessary conditions. Combining all of them will yield the statement. 

 We recall that for every $\xi\in k^*$  
\begin{equation}\label{equaz}
x_\alpha(\xi)x_{-\alpha}(-\xi^{-1})x_\alpha(\xi)\in  s_\alpha T.
\end{equation}

For $h\in k$ we consider the family of representatives of $C$ given by $y(h):=\theta^{-1}(\s_\alpha x_{\alpha}(h))*x$. Then, for
some structure constants $\eta_1,\eta_2,\eta_3, d_{\alpha\beta}$ that are always non-zero in good characteristic, and for some $t\in T$ we have:
\begin{align*}
y(h)&\in t \m \s_\alpha  x_{-\alpha}(\eta_1
h)x_{\alpha}(c_\alpha-h)x_\alpha(h)x_{\beta}(c_\beta)x_{\alpha+\beta}(c_{\alpha+\beta})x_{\alpha}(-h)\s_\alpha^{-1}P^u\\
&=t\m  x_\alpha(\eta_2
h)x_{-\alpha}(\eta_3(c_\alpha-h))\s_\alpha
x_{\beta}(c_\beta)x_{\alpha+\beta}(c_{\alpha+\beta}+hc_\beta
d_{\alpha\beta})\s_\alpha^{-1}P^u.
\end{align*}
Let $h_1$ and $h_2$ be the solutions of
\begin{equation*}X^2(\eta_2\eta_3)-c_\alpha
\eta_2\eta_3X-1=0\end{equation*} so that $-(\eta_2
h_i)^{-1}=(c_\alpha-h_i)\eta_3$ and we may apply \eqref{equaz}.
The elements corresponding to $h_1$ and $h_2$ satisfy
$$y(h_i)\in C\cap \m t' \s_\alpha x_{\beta}(\eta_4(c_{\alpha+\beta}+h_ic_\beta
d_{\alpha\beta}))P_\beta^u\subset C\cap B w_Cs_\alpha B$$
for some $t'\in T$ and some nonzero structure constant $\eta_4$. Here, $P_\beta$ denotes the minimal parabolic subgroup of $G_R$ associated with $\beta$.

We let now $\theta^{-1}(\s_\beta)$ act on $y(h_i)$ for $i=1,2$. We have, for some non-zero $\eta_5\in k$:
$$\theta^{-1}(\s_\beta)* y(h_i)\in B w_Cs_\beta s_{\alpha} s_\beta x_{-\beta}(\eta_5(c_{\alpha+\beta}+h_ic_\beta
d_{\alpha\beta})) B.$$ Moreover,  $w_Cs_\beta s_\alpha s_\beta\beta=\theta\alpha$ holds because $\alpha,\,\beta\in R_{w_C}$. Therefore, if we had $c_{\alpha+\beta}+h_ic_\beta
d_{\alpha\beta}\neq0$ we would have $\theta^{-1}(\s_\beta)* y(h_i)\in C\cap B w_Cs_\beta s_{\alpha}B$, contradicting the assumption on $C$ to be involutive. 
Thus
\begin{equation}\label{condiz}c_{\alpha+\beta}+h_ic_\beta
d_{\alpha\beta}=0.\end{equation} 
This condition must hold  for both $i=1,2$ thus we have either
$h_1=h_2$ so that
\begin{eqnarray}\label{Ialpha}&\Delta_\alpha=\eta^2_2\eta^2_3c_\alpha^2+4\eta_2\eta_3=0,\quad\mbox{
or }\\
\label{IIalpha}&c_\beta=c_{\alpha+\beta}=0.\end{eqnarray}

Let us now interchange the roles of $\alpha$ and $\beta$. We consider, for $l\in k$, the family of elements
\begin{align*}
z(l)&=\theta^{-1}(\s_\beta x_\beta(l))x x_{\beta}(-l)\s_\beta^{-1}\\
&\in \theta^{-1}(\s_\beta
x_\beta(l))\m
x_\beta(c_\beta)x_\alpha(c_\alpha)x_{\alpha+\beta}(c_{\alpha+\beta}+
c_\alpha c_\beta d_{\alpha\beta}) x_{\beta}(-l)\s_\beta^{-1}P^u.
\end{align*}

Using the same procedure as above with  $\beta$ and $\alpha$ interchanged we see that there are nonzero structure constants $\xi_1,\,\xi_2$, such that  if
$l_1$ and $l_2$  are the solutions of
\begin{equation*}\xi_1 X^2-c_\beta\xi_1 X-1=0\end{equation*} then
$$z(l_j)\in C\cap\m \s_\beta T x_\alpha(\xi_2(c_{\alpha+\beta}+c_\alpha
c_\beta d_{\alpha\beta}-l_jc_\alpha
d_{\alpha\beta}))P_\alpha^u,$$
for $j=1,\,2$, where $P^u_\alpha$ is as usual and $d_{\alpha\beta}$ is the structure constant occurring in \eqref{condiz}.

The action of $\theta^{-1}(\s_\alpha)$ on $z(l_j)$ for $j=1,2$ would yield an element in $C\cap Bw_C
s_\alpha s_\beta B$ unless
\begin{equation}\label{alpha+beta}c_{\alpha+\beta}+c_\alpha
c_\beta d_{\alpha\beta}-l_jc_\alpha
d_{\alpha\beta}=0\end{equation}
for both $j=1,2$. This forces either $l_1=l_2$ and therefore
\begin{eqnarray}\label{Ibeta}&\Delta_\beta=\xi^2_1c_\beta^2+4\xi_1=0,\quad\mbox{
    or }\\
\label{IIbeta}&c_\alpha=c_{\alpha+\beta}=0.
\end{eqnarray}
If $c_\alpha=0$ then \eqref{Ialpha} does not hold so $c_\alpha=c_\beta=c_{\alpha+\beta}=0$. In this situation any $k_{\alpha,\beta}$ will do.

If $c_\alpha\neq0$ then \eqref{Ibeta} must hold so we have at most two choices for $c_\beta$, and $c_\beta\neq0$. Thus, \eqref{Ialpha} must hold and we have a finite number of 
possibilities for $c_\alpha$, too. In this case, by \eqref{condiz}, we have  $c_{\alpha+\beta}=-\frac{1}{2}c_\alpha c_\beta d_{\alpha\beta}$ so we may take $k_{\alpha,\beta}=-\frac{1}{2}d_{\alpha\beta}$.\hfill$\Box$

\bigskip

\begin{lemma}\label{determined}Let $C$ be an involutive $\theta$-twisted conjugacy class and let $\m\in N(T)$ be such that $C\cap \m U\neq\emptyset$.
Let $\Delta_R=\{\gamma_1,\ldots,\,\gamma_r\}$ and let $x=\m\v\in C\cap \m U$.
Then, for every $\gamma=\sum_{j=1}^r n_j\gamma_j\in R_{w_C}$  there is a polynomial $p_{\m,\gamma}(X)\in k[X_j~|~n_j\neq0]$ without constant term, 
depending only on  $\gamma$, $\m$,  the fixed ordering of the positive roots in $R_{w_C}$, and the structure constants of $G$, 
such that the coefficient $c_\gamma$ of $x_\gamma$ in the expression of $\v$ is the evaluation of $p_{\m,\gamma}(X)$ at $X_j=c_{\gamma_j}$ for every $j=1,\,\ldots,\,r$ in the support of $\gamma$.  In particular,  we have $|C\cap mU|<\infty$ for every $m\in w_CT$.
\end{lemma}
\pf  Without loss of generality we may assume that the ordering of the positive roots is with increasing height. We shall proceed by induction on the height $h$ of the root $\gamma$ with respect to $\Delta_R$.
Let us assume that the claim holds for all $\gamma$ with ${\rm ht}\,\gamma\leq h-1$. Let $\nu\in R_{w_C}$ with ${\rm
  ht}\,\nu=h$. By Lemma \ref{alpha-beta}  the statement holds for $h\leq 2$, so we will assume that $h$ is greater than $2$. There exists $\beta\in\Delta_R$ for which
${\rm ht}\,s_\beta\nu= h-1$. The strategy will be to consider $y=\theta^{-1}(\s_\beta)*x$ and to find an element $z$ in $X_{\theta^{-1}(\beta)}*y$ lying in $w_CB$. The induction step will be obtained by comparing the coefficient of $x_\nu$ in the expression of $\v$ with the coefficient of $x_{s_\beta(\nu)}$ in the expression of $u$, for $z=\m tu$. We have 
\begin{equation}y=(\theta^{-1}\s_\beta)* x=(\theta^{-1}\s_\beta) \m \s_\beta^{-1} (\s_\beta\v   \s_\beta^{-1})=\m t
\prod_{\gamma\in R_{w_C}}x_{s_\beta\gamma}(\eta_\gamma c_\gamma)\end{equation}
 for some non-zero structure constants $\eta_\gamma$ and some $t\in T$ depending on $\s_\beta$ and $\m$. Here, the product
respects  the fixed ordering of the $\gamma$'s and not of the $s_\beta\gamma$'s. We have:
$y=\m t \v_1 x_{-\beta}(\eta_\beta c_\beta)\v_2$ for some $\v_1,\,\v_2\in P^u$, the unipotent radical of the minimal standard parabolic subgroup $P$ of $G_R$ associated with $\beta\in \Delta_R$. Since the ordering is increasing in height, $\v_1$ is a product of elements of the form $x_{s_\beta\gamma}(\eta_\gamma c_\gamma)$ for $\gamma\in\Delta_R$.

Let $\eta\in k$ be  such that $x_{\theta^{-1}\beta}(\eta c_\beta)\m t=\m t x_{-\beta}(-\eta_\beta c_\beta)$ and let us consider the element $z=x_{\theta^{-1}\beta}(\eta c_\beta)*y$. Then
$$\begin{array}{rl}z&=\m t(x_{-\beta}(-\eta_\beta c_\beta)\v_1 x_{-\beta}(\eta_\beta c_\beta))\v_2 x_\beta(-\eta c_\beta)=\m t u\\
&=\m t \prod_{\gamma\in R_{w_C}}x_\gamma(d_\gamma)\in\m t U\cap C
\end{array}$$ where the product is taken according to the ordering of the positive roots in $R_{w_C}$. Here, the expression $x_{-\beta}(-\eta_\beta c_\beta)\v_1 x_{-\beta}(\eta_\beta c_\beta)$ is a product of terms of the form $x_{s_\beta\gamma}(\eta_\gamma c_\gamma)$ for $\gamma\in\Delta_R$ such that $s_\beta\gamma=\gamma$ and  terms of the form
$x_{s_\beta\gamma'}(\eta_{\gamma'} c_{\gamma'})x_{s_\beta\gamma'-\beta}(\eta_{\gamma',\beta} c_{\gamma'} c_\beta)$ with $\eta_{\gamma',\beta}$ a nonzero structure constant for $\gamma'\in\Delta_R$ such that $s_\beta\gamma'=\gamma'+\beta$.

By the induction hypothesis applied to $z$ and $s_\beta\nu$, the coefficient $d_{s_\beta\nu}$ is evaluation at the 
 $d_\alpha$ for $\alpha$ in the support of $s_\beta\nu$  of a polynomial $p_{\m t,s_\beta\nu}(X)$ without constant term.
Besides, each
$d_{\mu}$ differs from $\eta_{s_\beta\mu} c_{s_\beta\mu}$ by a
 (possibly trivial) sum of monomials in the $c_{\mu'}$, $c_\beta$, the structure constants $\eta_{\mu'}$, $\eta_\beta$, and the structure constants coming from application of
 Chevalley's formula \cite[Proposition 8.2.3]{springer}  when reordering root subgroups. More precisely, we have 
\begin{equation}\label{reordering}d_\mu=\eta_{s_\beta\mu} c_{s_\beta\mu}+\sum C_{i_1,\ldots, i_p,j}^{j_1,\ldots,j_q}(\prod_{l=1}^pc_{\nu_l}^{i_l})c_\beta^j\prod_{r=1}^q(c_{\gamma'_r}c_\beta)^{j_r}
\end{equation}
where $C_{i_1,\ldots, i_p,j}^{j_1,\ldots,j_q}$ denotes a coefficient depending on the structure constants. The sum is taken over the possible decompositions 
$$\mu=\sum_{l=1}^pi_ls_\beta\nu_l+j\beta+\sum_{r=1}^qj_r(s_\beta\gamma'_r-\beta)$$ for $i_l>0,\,j_r>0$ and $j\geq0$
and $\gamma'_r\in\Delta_R$ such that $s_\beta\gamma'_r=\gamma'_r+\beta$. Contribution to $d_{s_\beta\nu}$ as in \eqref{reordering} may occur only when
 \begin{equation}\label{nu}s_\beta\nu=\sum_{l=1}^pi_ls_\beta\nu_l+j\beta+\sum_{r=1}^qj_r(s_\beta\gamma'_r-\beta) \end{equation} for $i_l, j_r>0$ and
 $j\geq0$. Then, ${\rm ht}s_\beta\nu_l<{\rm ht}s_\beta\nu= h-1$ and $s_\beta\gamma'_r-\beta=\gamma'_r$. 
 Applying $s_\beta$ to \eqref{nu} we have
 \begin{equation}\label{nu-s}\nu+j\beta=\sum_{l=1}^pi_l\nu_l+\sum_{r=1}^qj_r(\gamma'_r+\beta) \end{equation}
 so the support of $\nu$ contains $\gamma'_r$ and the support of $\nu_l$.  Since $\Phi$ is simply-laced, ${\rm ht}\nu_l\leq  {\rm ht}s_\beta\nu_l+1<h$ for every $l$. 
We may thus  apply the induction hypothesis to $c_{\nu_l}$. So
$$c_\nu=\eta_\nu^{-1}p_{\m t, s_\beta\nu}(d_{\gamma_i})-\eta_\nu^{-1}\sum C_{i_1,\ldots, i_p,j}^{j_1,\ldots,j_q}\left(\prod_{l=1}^p(p_{\m ,\nu_l}(c_{\gamma_i}))^{i_l}\right)c_\beta^j\prod_{r=1}^q(c_{\gamma'_r}c_\beta)^{j_r}.$$
The statement is proved if we show that $d_{\gamma_i}$ is a monomial in  $c_{\gamma_i}$ and possibly $c_\beta$ multiplied by a structure constant. 
If $s_\beta(\gamma_i)=\gamma_i$  then the coefficient $d_{\gamma_i}$ of  $x_{\gamma_i}$ in $u$ is equal to $\eta_{\gamma_i}c_{\gamma_i}$. 
If, instead, $s_\beta(\gamma_i)=\gamma_i+\beta$ and $\gamma_i$ follows $\beta$ in the ordering of the positive roots, then the coefficient $d_{\gamma_i}$ of $x_{\gamma_i}$ in $u$ equals $\eta_{\gamma_i+\beta}c_{\gamma_i+\beta}$, which is as required by Lemma \ref{alpha-beta}.
Finally, if $s_\beta(\gamma_i)=\gamma_i+\beta$ and $\gamma_i$  precedes $\beta$ in the ordering of the positive roots, then the coefficient $d_{\gamma_i}$ of $x_{\gamma_i}$ in $u$ is the sum of  $\eta_{\gamma_i+\beta}c_{\gamma_i+\beta}$ with the correction term obtained from moving $x_{-\beta}(\eta_\beta c_\beta)$ from the right of $\v_1$ to the left. 
By Chevalley's commutator formula,  this correction term equals $Kc_{\gamma_i}c_\beta$ for some product $K$ of structure constants. 
Thus, $c_\nu$ is evaluation of a polynomial without constant term depending only on the structure constants, on the choice of $\m$, and the fixed ordering of the roots.
The last statement follows from Lemma \ref{alpha-beta} because when $w_0\neq-1$, if $R$ is non-empty, then it is always irreducible of rank greater than $1$.
\hfill$\Box$
\smallskip

\begin{remark}{\rm It follows from Lemma \ref{determined} that  if $\v\in P_\alpha^u$ for
{\em every} $\alpha\in\Delta_R$ then $\v=1$ so $x=\m$.  In particular,  by the proof of Lemma \ref{alpha-beta},  this condition holds if  $\v\in P_\alpha^u$  for {\em some} $\alpha\in \Delta_R$ because $R$ is either.
empty or irreducible of rank greater than $1$.}
\end{remark}

Combining Lemma \ref{sufficient}, Proposition \ref{meets},  and Lemma \ref{determined} we have the following result.

\begin{theorem}\label{quasitutti}Let $C$ be an involutive $\theta$-twisted conjugacy class in $G$. If for $w_C=w_0w_\Pi$ we have $(\Phi,\Pi)\neq(D_{2n},\emptyset)$, then $C$ is spherical.
\end{theorem}


\section{The case of $(D_{2n},\emptyset)$}\label{remaining}

Let us now consider the case of involutive $\theta$-twisted conjugacy classes $C$ in type $D_{2n}$ with $w_C=w_0$. In this case 
$\Delta_R=\{\alpha_1,\,\ldots,\,\alpha_{2n-2},\,\alpha_{2n-2}+\alpha_{2n-1}+\alpha_{2n}\}$ whereas $C_{w_0}\cap\Delta=\{\alpha_{2n-1},\,\alpha_{2n}\}$. Let $G_R'=[G_R,\,G_R]$ where $G_R$ is as in Section \ref{involutive}. The automorphism $\theta$ stabilizes $G_R'$ and it acts trivially on it. 

\begin{lemma}\label{piccola}Let $C$ be an involutive $\theta$-twisted conjugacy class in type $D_{2n}$ with $w_C=w_0$ and let $x=\w_0^{-1}\v\in T w_0B\cap C$. Then, the $G_R'$-orbit $G_R'\cdot_\theta x$ is spherical.
\end{lemma}
\pf By Remark \ref{slodowy} the semisimple group
$G_R'$ is simply-connected. It is in fact simple of type $D_{2n-1}$. 
Let $B_R=B\cap G'_R$ and let $T_R$ be the maximal torus of $G_R'$, generated by the elements of the form $h_{\gamma}(\xi)$ for $\gamma\in\Delta_R$ (see Section \ref{notation}). 
By Proposition \ref{real} we have  $x=\w_0^{-1}\v\in Tw_0U_R\cap C$ . We may choose a representative $\w_R$ of the longest element of the Weyl group of $G_R$ in $N(T)\cap G_R'$.
Conjugation by $\w_R\w_0$ stabilizes $T$, $B_R$ and $T_R$ and it induces a non-trivial automorphism of $R$. Thus, for some $t\in T_R$, conjugation by $t\w_R\w_0$ is the automorphism induced by the non-trivial involution $\tau$ of the Dynkin diagram of $G'_R$. Let $g\in G_R'$. We have  $\theta(g)=g$ so $g\cdot_\theta x =g x g^{-1}$  and the morphism
$$\begin{array}{rl}
f\colon G'_R\cdot_\theta x&\to G'_R\cdot_\tau (t\w_R\v)^{-1}\\
z&\mapsto (t\w_R\w_0 z)^{-1}
\end{array}
$$
is  a $G'_R$-equivariant isomorphism. So, it is enough to show that the $G_R'$-variety $G'_R\cdot_\tau (t\w_R\v)^{-1}$ is spherical. 
We shall show that $G_R'\cdot_\tau (t\w_R\v)^{-1}$ is involutive. The statement will follow from Theorem \ref{quasitutti}.
Let $\sigma$ be an element in the Weyl group of $G_R'$ such that $G_R'\cdot_\tau (t\w_R\v)^{-1}\cap B_R\sigma B_R\neq\emptyset$ and
let $y\in G_R'\cdot_\tau (t\w_R\v)^{-1}\cap B_R\sigma B_R$. Then 
$$\begin{array}{rl}
f^{-1}(y)&=\w_0^{-1}\w_R^{-1}t^{-1}y^{-1}\in \w_0^{-1}\w_R^{-1} t^{-1}B_R\sigma^{-1} B_R\\
&= B_R \w_0^{-1}\w_R^{-1} \sigma^{-1} B_R\subset B w_0^{-1}w_R^{-1} \sigma^{-1}B\cap C.
\end{array}$$
Since $C$ is $\theta$-involutive and $w_R$ and  $\sigma$ are all $\theta$-invariant because they are products of $\theta$-invariant reflections, we have
$w_0^{-1}w_R^{-1} \sigma^{-1}w_0^{-1}w_R^{-1} \sigma^{-1}=1$. The involutions $w_0$ and $w_R$ commute because $w_0$ acts trivially on each reflection in $W_{\Delta_R}$ so 
$(w_Rw_0\sigma^{-1}w_0^{-1} w_R^{-1})\sigma^{-1}=1$ and  $\sigma\tau(\sigma)=1$. 
\hfill$\Box$

\bigskip

\begin{lemma}\label{real-in}Let $C$ be an involutive $\theta$-twisted conjugacy class in type $D_{2n}$ with $w_C=w_0$ and let $x=\w_0\v\in Tw_0B\cap C$. Let $\alpha\in\Delta_R$. Then for all but finitely many $\xi\in k$ the set $x_\alpha(\xi)s_\alpha B\cap G_x$ is non-empty.
\end{lemma}
\pf  We have $\theta\alpha=\alpha$ and $\theta(x_\alpha(\xi))=x_\alpha(\xi)$. Let  $\v=x_\alpha(c)\v'\in x_\alpha(c)P_\alpha^u$ and let us consider the following representatives of $C$, for $\xi\in k$: 
$$y_\xi=x_{\alpha}(\xi)\cdot_\theta x=x_\alpha(\xi)\w_0\v x_{\alpha}(-\xi)=\w_0 x_{-\alpha}(\eta\xi)x_{\alpha}(c-\xi)\v'' $$ for some nonzero structure constant $\eta$ and some $\v''\in P_\alpha^u$,
and 
$$\begin{array}{rl}
z_{\xi}&=\s_\alpha\cdot_\theta y_\xi=\s_\alpha \w_0 x_{-\alpha}(\eta\xi)x_\alpha(c-\xi)\v'' \s_\alpha^{-1}\\
&\in \w_0 t x_{\alpha}(\eta'\xi)x_{-\alpha}(\eta''(c-\xi))P_\alpha^u
\end{array}$$
for some nonzero structure constants $\eta',\eta''$ and some $t\in T$.
It follows from \eqref{equaz} that if $\eta'\xi\eta''(c-\xi)\neq-1$ then $z_\xi\in Bw_0B$.

Let $v_0$ be the dense $B_R$-orbit in $G'_R\cdot_\theta x$ and let $B\sigma B$ be the Bruhat double coset in $G$ containing it. Then 
$\overline{G_R'\cdot_\theta x}=\overline{v_0}\subset (\bigcup_{\omega\leq\sigma}B\omega B)$. Since $Bw_0B\cap G_R'\cdot_\theta x$ is non-empty we necessarily have $\sigma=w_0$. 

Moreover, using uniqueness of the Bruhat decomposition as in \cite[Lemma 2.1]{lu} , \cite[Theorem 4.1]{mauro} or \cite[Theorem 5]{ccc},  we may show that the $\theta$-centralizer in $B_R\subset U_RT$ of an element in $w_0B$ is finite. This shows that every $B_R$-orbit $v$ in  $G_R'\cdot_\theta x$ which is contained in $Bw_0B$ has the same dimension  as the dense one, that is $\dim v=\dim B_R=\dim v_0$. Therefore, $v$ must coincide with $v_0$.  Thus, $z_{\xi}$ and $x$ lie in $v_0$ and there is $b_\xi\in B_R$ such that $b_\xi\cdot_\theta z_{\xi}=x$. In other words, for every $\xi$ but finitely many there is an element in 
$G_x\cap B_R \s_\alpha x_{\alpha}(\xi)\subset G_x\cap B \s_\alpha x_\alpha(\xi)$. Taking inverses we have the statement.\hfill$\Box$

\begin{lemma}\label{complex}Let $C$ be an involutive $\theta$-twisted conjugacy class in type $D_{2n}$ with $w_C=w_0$ and let $x=\w_0\v \in w_0TU\cap C$. Let $\alpha\in\Delta\cap C_{w_0}$. Then for every $\xi\in k$ the set $x_{-\alpha}(\xi)U\cap G_x\neq\emptyset$.
\end{lemma}
\pf The simple root $\alpha$ is either $\alpha_{2n-1}$ or $\alpha_{2n}$ so $\alpha\pm\theta\alpha\not\in\Phi$. Let us consider the following representatives of $C$ for $\xi\in k$:
$$y_\xi=x_{\theta\alpha}(\xi)x x_{\alpha}(-\xi)=\w_0 x_{-\theta\alpha}(\eta\xi)\v x_{\alpha}(-\xi)$$ for some non-zero structure constant $\eta$ and
$$z_\xi=x_{-\alpha}(\eta\xi)y_\xi x_{-\theta\alpha}(-\eta\xi)=\w_0 x_{\alpha}(\eta'\xi)x_{-\theta\alpha}(\eta\xi)\v x_{-\theta\alpha}(-\eta\xi)  x_{\alpha}(-\xi)$$
for some nonzero structure constant $\eta'$.

The element $x_{\alpha}(\eta'\xi)x_{-\theta\alpha}(\eta\xi)\v x_{-\theta\alpha}(-\eta\xi)  x_{\alpha}(-\xi)$ lies in $U$ because 
$\v$ lies in $P_{\theta\alpha}^u$ by Proposition \ref{real}. Applying the Proposition once more to $z_\xi\in C\cap \w_0U$ we see that $\eta'=1$. 

Let us fix an ordering of the positive roots so that all $\theta$-invariant roots precede the non-invariant ones. We recall that the non-invariant positive roots are of the form $\beta_i=\alpha_i+\cdots+\alpha_{2n-2}+\alpha$ or   $\theta\beta_i$. It follows from Chevalley's commutator formula \cite[Proposition 8.2.3]{springer} that $x_{-\theta\alpha}(\eta\xi)\v x_{-\theta\alpha}(-\eta\xi)=\v\v'$ where $\v'$ lies  in the abelian subgroup $Y_\alpha$ of $U$ generated by the root subgroups associated with the $\beta_i$'s. Besides, $x_{\alpha}(\xi)\v \v'  x_{\alpha}(-\xi)=\v\v'\v''$ where $\v''$ lies again in $Y_\alpha$.
By Proposition \ref{real}  we conclude that $\v'\v''=1$ so $z_\xi=x$ and  $x_{-\alpha}(\eta\xi)  x_{\theta\alpha}(\xi)\in G_x$. Since $\eta$ is a fixed non-zero structure constant and the statement holds for every $\xi$, we have the statement.
\hfill$\Box$

\begin{lemma}\label{dense}Let $X$ be a transitive $G$-variety and let $x\in X$. If, for every $\alpha\in \Delta$ we have $x_\alpha(\xi) s_\alpha B\cap G_x\neq\emptyset$ for all but finitely many $\xi\in k$, then the space $X$ is spherical and $B\cdot x$ is the dense $B$-orbit.
\end{lemma}
\pf It is enough to show that $BG_x$ or, alternatively, $G_xB$, is dense in $G$. We will do so by showing that $G_xB\cap B w_0 B$ is dense in $Bw_0B$. 

Let $U^w$ be the subgroup generated by the root subgroups associated with roots in $\Phi_w=\{\alpha\in \Phi^+~|~w^{-1}\alpha\in-\Phi^+\}$. 
Then $BwB=U^wwB$ and, once we have fixed an ordering of the roots in $\Phi_w$,  we may identify $U^w wB/B\subset G/B$ with 
 the affine space ${\mathbb A}^{\ell(w)}$ through
the map $u w B=\prod_{\gamma\in\Phi_w}x_\gamma(c_\gamma)w B\mapsto
(c_\gamma)_{\gamma\in\Phi_w}$. We will show by induction on the length $\ell(w)$ of $w$ that  the set $U^w_0$ of elements $u$ in $U^w$ for which
$u w B\cap G_x$ is non-empty contains the complement of the union of finitely many hyperplanes in $U^w\cong {\mathbb A}^{\ell(w)}$.
For $w=1$ there is nothing to say. Suppose that
the claim holds for $\ell(w)=l$. We consider
$\omega\in W$ with $\ell(\omega)=l+1$. Then $\omega=\sigma s_\alpha$ for
some $\sigma\in W$ with $\ell(\sigma)=l$ and some $\alpha\in\Delta$ with $\sigma\alpha\in\Phi^+$. 
Besides, $\Phi_\omega=\Phi_\sigma\cup\{\sigma\alpha\}$ so $U^\omega=U^\sigma X_{\sigma\alpha}$. 

By the hypothesis, for all but finitely many $\xi\in k$ and  for every $u\in U_0^\sigma$ there is a $b\in B$ for which $(u\sigmad b)(x_{\alpha}(\xi)\s_\alpha B)\cap G_x\neq\emptyset$. Let $b=x_{\alpha}(r)\v$ for $r\in k$ and $\v\in P^u_\alpha$. Then for some nonzero structure constant $\eta$ we have
$$(u\sigmad b)(x_{\alpha}(\xi)s_\alpha B)=u \sigmad x_{\alpha}(r+\xi)s_\alpha B=u x_{\sigma\alpha}(\eta(r+\xi))\sigma s_\alpha B$$ and 
$u x_{\sigma\alpha}(\eta(r+\xi))\sigma s_\alpha B\cap G_x\neq\emptyset$.
Since all but finitely many $\xi$ were allowed and $\eta\neq0$  the intersection $G_x\cap B\omega B$ contains $U^\sigma_0 x_{\sigma\alpha}(\xi)\omega B$ for all but finitely many $\xi$, thus $U_0^\omega$ contains the complement of finitely many hyperplanes in ${\mathbb A}^{\ell(\omega)}$. \hfill$\Box$

\smallskip

\begin{proposition}\label{casoaparte}Let $C$ be an involutive $\theta$-twisted conjugacy class in type $D_{2n}$ with $w_C=w_0$ . Then $C$ is spherical.
\end{proposition}
\pf We have $\Delta=(\Delta\cap R_{w_0})\cup(\Delta\cap C_{w_0})$. By Lemmas \ref{real-in}, \ref{complex} and formula \eqref{equaz} the hypotheses of Lemma \ref{dense} are satisfied.\hfill$\Box$

\bigskip

\begin{remark}\label{general}{\rm Let $\pi\colon G\to H$ be a central isogeny of simple groups with $G$ simply-connected and suppose that the automorphism $\theta$ of $G$ preserves ${\rm Ker}(\pi)$. Then it preserves the character group of the maximal torus $T_H=\pi(T)$ of $H$. The automorphism $\theta$ induces an automorphism $\overline{\theta}$ of $H$ such that $\overline{\theta}\circ\pi=\pi\circ\theta$. 
Thus, the $\theta$-twisted conjugacy classes of $G$ are mapped onto $\overline{\theta}$-twisted conjugacy classes of $H$. Clearly, the Bruhat cells they intersect correspond to the same Weyl group elements. Moreover, 
 it is not hard to verify that 
%
$C$ is spherical if and only if $\pi(C)$ is so. This allows the generalization of the obtained results from simply-connected groups to  a more general setting.}
\end{remark}

\bigskip

Combining Theorem \ref{induction}, Theorem \ref{quasitutti}, Proposition \ref{casoaparte} and Remark \ref{general} we obtain our main result.

\begin{theorem}Let $G$ be a simple algebraic group over an algebraically closed field of good odd characteristic. Let $B$ be a Borel subgroup of $G$ and $T$ a maximal torus in $B$. Let $\theta$ be an involution of the Dynkin diagram of $G$ preserving the character group of $T$. A $\theta$-twisted conjugacy class 
$C$ is spherical if and only if it lies in $\bigcup_{w\theta(w)=1}BwB$. 
\end{theorem}

\begin{remark}{\rm The above theorem can be viewed as an analogue of \cite[Theorem 5.7]{gio-fourier} for non-connected semisimple groups with simple identity component}. 
\end{remark}

\section{The dimension formula in good characteristic}\label{dimension-section}

In this section we will show how we get \cite[Theorem 1.1]{lu} in good, odd characteristic, for $\theta$ a non-trivial involution of the Dynkin diagram, as a by-product of the  results in the previous sections.

\begin{proposition}\label{dimension-maximal}Let $B* x$ be a maximal $B$-orbit in an involutive $\theta$-twisted conjugacy class $C$. Then 
\begin{equation}\dim B* x=\ell(w_C)+{\rm rk}(1- w_C\theta).\end{equation} 
\end{proposition}
\pf Let  us choose $x=\m\v\in C\cap  Tw_CU$. We recall that $\Phi^+_\Pi$ is the set of roots whose positivity is preserved by $w_C$ so $\dim U_\Pi=|\Phi^+|-\ell(w_C)$.
It follows from uniqueness of the Bruhat decomposition (see \cite[Theorem 4.1]{mauro} or \cite[Lemma 2.1]{lu}) that the $\theta$-stabilizer of $x$ in $B=U^{w_C}U_\Pi T$, with notation as in the proof of Lemma \ref{dense}, is contained in $U_\Pi T^{w_C\theta}$. Hence,  $\dim B* x\geq \ell(w_C)+{\rm rk}(1-w_C\theta)$. 

Let $\alpha\in\Pi$. Then
$$x_\alpha(t)* \m \v= (x_{\alpha}(t)*\m)(x_{\theta\alpha}(t)\v x_{\theta \alpha}(-t)).$$
By Lemma \ref{levi-commute} we have $x_{\alpha}(t)*\m=\m$ and by Proposition \ref{real} we have $x_{\theta\alpha}(t)\v x_{\theta\alpha}(-t)=\v$, so $x_\alpha(t)$ lies in the $\theta$-stabilizer of $x$, and therefore the same holds for all elements in $U_\Pi$. Moreover, the maximal torus $T_\Pi$ of $[L_\Pi,L_\Pi]$ generated by the $h_\alpha(\zeta)$ for $\alpha\in\Pi$ is contained in $(T^{w_C\theta})^\circ$. It is not hard to verify by a dimensional argument that $(T^{w_C\theta})^\circ$ is equal to $T_\Pi$ for all choices for $w_C$ except from $(D_{2n},\emptyset)$. Since $\Pi\perp R$, Lemma \ref{levi-commute}  and Proposition \ref{real} imply
$B_x^\circ=(T^{w_C\theta})^\circ U_\Pi$ and the statement. 

Let us now consider the case $(D_{2n},\emptyset)$. 
In this  case $B_x\subset T^{w_C\theta}$ and $(T^{w_C\theta})^\circ$ is the $1$-dimensional torus of the elements $h_\xi=h_{\alpha_{2n-1}}(\xi)h_{\alpha_{2n}}(\xi^{-1}),$ for $\xi\in k^*$. These elements certainly lie in the $\theta$-stabilizer of $\m$. 
We have $h_\xi* x=(h_\xi*\m) \theta(h_\xi)\v \theta(h_\xi)^{-1}=\m h_\xi^{-1}\v h_\xi$. By Proposition \ref{real} the element $h_\xi$ centralizes $\v$ because the roots in $R_{w_0}$ are orthogonal to the $-1$ eigenspace of $\theta$.  We have $(B_x)^\circ=(T^{w_C\theta})^\circ$ and the statement.
\hfill$\Box$

\bigskip
The main result of this section follows:
\bigskip
\begin{theorem}\label{dimension-formula}Let $G$ be a simple group over an algebraically closed field of good, odd characteristic. Let $\theta$ be an involution of its Dynkin diagram and let us assume that the character group of $T$ is $\theta$-invariant. Then, a $\theta$-twisted conjugacy class $C$ is spherical if and only if $\dim C=\ell(w_C)+{\rm rk}(1-w_C\theta)$.
\end{theorem}
\pf Let us assume first that $G$ is simply-connected. If $C$ is spherical then its dense $B$-orbit $v_0$ is necessarily maximal so $\dim C=\dim v_0$. Moreover, $C$ is involutive by Theorem \ref{induction}  so Proposition \ref{dimension-maximal} yields the statement in this case. For the general case we use Remark \ref{general}.

If $\dim C=\ell(w_C)+{\rm rk}(1-w_C\theta)$ we argue  as in \cite{ccc},\cite{mauro} or \cite{lu}. The idea is to take a representative $x$ of $C$ in $w_CB$ and use uniqueness of the Bruhat decomposition in order to show that  $B_x$ is contained in $T^{w_C\theta}U_\Pi$. This implies that the dimension of the $B$-orbit of $x$ equals the dimension of $C$.\hfill$\Box$

\begin{remark}{\rm The dimension formula in \cite{lu} is stated in characteristic zero and it generalizes to $\theta$ non-trivial the dimension formula in \cite{ccc,gio-mathZ}. The proof works also in positive characteristic provided that some requirements on the base field listed in \cite[Remark 2.3]{lu}  hold. The present proof covers also the case in which the characteristic of $k$ is not {\em very good}, i.e., ${\rm char} k$ divides $n+1$ in type $A_n$. In this case, the orbit map to a twisted conjugacy class is not necessarily separable (\cite[Page 380]{mohr}), so the requirements in \cite[Remark 2.3]{lu} are not satisfied. On the other hand, \cite{lu} covers the triality case, whereas the present approach does not reach the case of triality with $w_C=w_0s_2$. The dimension formula in the triality case when $w_C=w_0$ easily follows from the fact that $(T^{w_C\theta})^\circ=1$ in this case, so the argument in \cite[Theorem 5]{ccc} already shows that the dimension of a $B$-orbit in $Bw_0B$ is equal to the dimension of $B$, which is equal to $\ell(w_0)+{\rm rk}(1-w_0\theta)$. }
\end{remark}

\section{Acknowledgements}
I wish to thank Kei Yuen Chan for pointing out an inaccuracy in a first version of the paper, Jiang-Hua Lu for answering several questions, and the referees for careful reading, pointing out critical points in the proof of Lemma~\ref{determined}   and
many other useful suggestions, corrections and remarks.  This paper was typed on a computer supported by  Project CPDA071244 of the University of Padova.

\end{document}